# Asymptotic behavior "almost everywhere" of additive and multiplicative arithmetic functions

Victor Volfson

ABSTRACT  We define the asymptotic behavior "almost everywhere" of additive and multiplicative arithmetic functions in the paper. Classes of additive and multiplicative arithmetic functions are singled out for which asymptotics coincides "almost everywhere" with the asymptotics of the corresponding strongly additive and multiplicative arithmetic functions. Several assertions are proved and examples are considered.

## 1. INTRODUCTION

An arithmetic function is a function defined on the set of natural numbers and taking values on the set of complex numbers. The name - arithmetic function is due to the fact that this function expresses some arithmetic property of the natural series.

Estimating asymptotics of arithmetic functions and moments of arithmetic functions has been and remains an urgent problem at the present time [1], [2], [3].

Let's look at some properties of arithmetic functions.

Let a natural number have a canonical decomposition $m = p_1^{a_1}...p_t^{a_t}$, where $p_i$ is prime mumber, and $\alpha_i$ is natural number.

Then, the property is fulfilled for the additive arithmetic function:
$$f(m) = f(p_1^{a_1}...p_t^{a_t}) = f(p_1^{a_1}) + ... + f(p_t^{a_t}) = \sum_{p^\alpha \| m} f(p^\alpha).$$

The following property holds for the corresponding strongly additive function:
$$f^*(m) = f^*(p_1^{a_1}...p_t^{a_t}) = \sum_{p|m} f(p).$$

The property is executed for the multiplicative arithmetic function:
$$g(m) = g(p_1^{a_1}...p_t^{a_t}) = \prod_{p^\alpha \| m} g(p).$$

The following property holds for the corresponding strongly multiplicative function:
$$g^*(m) = g^*(p_1^{a_1}...p_t^{a_t}) = \prod_{p|m} g(p).$$

A multiplicative arithmetic function $g(m) > 0$ is converted to an additive one using the logarithm operation:
$$\ln(g(m)) = \ln(g(p_1^{a_1}...p_t^{a_t})) = \ln(\prod_{p^\alpha \| m} g(p)) = \sum_{p^\alpha \| m} \ln(g(p)).$$

An additive arithmetic function is converted into a multiplicative one using the potentiation operation : $e^{f(m)} = e^{f(p_1^{a_1}...p_t^{a_t})} = e^{f(p_1^{a_1})+...+f(p_t^{a_t})} = \prod_{p^\alpha \| m} e^{f(p^\alpha)}$.

Logarithm and potentiation can be carried out in any base $a > 0, a \neq 1$. We will consider $a = e$ in this work.

Arithmetic functions usually have the specific property of changing internally irregularly and chaotically, as a result of which the classical methods of analysis, as a rule, are powerless to adequately describe their behavior.

Probabilistic number theory responds to the problem naturally arising in such circumstances, to conduct a statistical study of this behavior.

Along with the concepts of extreme and average orders, which allow a cursory classification, there is the concept of the normal distribution of an arithmetic function, designed to reflect its "almost" certain behavior.

In practice, this involves the process of neglecting a set of zero-density integers (obviously depending on the function in question) to eliminate erroneous values.

Strikingly, this approach leads to the fact that order and regularity suddenly emerge from the apparent internal chaos. The illuminating focus of the "almost everywhere" concept opens up a new area of research that requires specific methods and yields specific results.

This approach to determining the asymptotic of an arithmetic function is considered in the paper.

## 2. ASYMPTOTIC EVALUTATION OF AN ARITHMETIC FUNCTION

Any initial segment of the natural series $\{1,2,...,n\}$ can be turned into a probability space $(\Omega_n, \mathcal{A}_n, \mathbb{P}_n)$ by taking $\Omega_n = \{1,2,...,n\}$, $\mathcal{A}_n$ — all subsets of $\Omega_n$, $\mathbb{P}_n(A) = \frac{|A|}{n}$. Then an arbitrary (real) function of the natural argument $f(m)$ (more precisely, its restriction to $\Omega_n$) can be considered as a random variable $\xi_n$ on this probability space: $\xi_n(m) = f(m)$, $1 \leqslant m \leqslant n$.

We can talk about the average value of the arithmetic function $f(m)$ on this probability space

$$E[f,n] = \frac{1}{n}\sum_{m=1}^{n} f(m), \qquad (2.1)$$

and dispersion of the arithmetic function :

$$D[f,n] = \frac{1}{n}\sum_{m=1}^{n} |f(m)|^2 - E^2[f,n]. \qquad (2.2)$$

Therefore, one can write Chebyshev's inequality for an arithmetic function $f(m), m = 1,...,n$ on a given probability space:

$$P_n(|f(m) - E[f,n]| \leq b\sqrt{D[f,n]}) \geq 1 - 1/b^2, \qquad (2.3)$$

where $b \geq 1$.

We put $b = b(n)$ in (2.3), where $b(n)$ is an indefinitely increasing function at value $n \to \infty$ and we get:

$$P_n(|f(m) - E[f,n]| \leq b(n)\sqrt{D[f,n]}) \geq 1 - 1/b^2(n). \qquad (2.4)$$

Limit $P_n$ (2.4) exists for value $n \to \infty$:

$$P_n(|f(m) - E[f,n]| \leq b(n)\sqrt{D[f,n]}) \to 1, n \to \infty. \qquad (2.5)$$

Expression (2.5) is an analogue of the law of large numbers for an arithmetic function [4].

Based on (2.5), we can write an expression for the asymptotic of the arithmetic function $f(m), m = 1,...,n$, which is satisfied almost everywhere for the value $n \to \infty$:

$$f(n) = E[f,n] + O(b(n)\sqrt{D[f,n]}), \qquad (2.6)$$

where $b(n)$ is a slowly growing function, i.e. $b(n) \leq n^\xi$, a $\xi > 0$.

Taking into account (2.6), if $b(n)\sqrt{D[f,n]} = o(E[f,n])$, then the following holds almost everywhere:

$$f(n) = E[f,n](1 + o(1)), \qquad (2.7)$$

otherwise

$$f(n) = O(b(n)\sqrt{D[f,n]}) \qquad (2.8)$$

at value $n \to \infty$.

Therefore, it is necessary to know the values $E[f,n]$ and $D[f,n]$ to determine the asymptotic of the arithmetic function $f(m), m = 1,...,n$ for the value $n \to \infty$.

Let's look at the asymptotic of the mean value for an additive arithmetic function $f(m), m = 1, ..., n$.

Assertion 1

The asymptotic of the mean value of an additive arithmetic function $f(m), m = 1, ..., n$ on the interval $[1, n]$ at value $n \to \infty$ is:

$$A_n = \sum_{p \leq n} \frac{f(p^\alpha)}{p^\alpha}. \qquad (2.9)$$

Proof

We define a random variable $f^p(m) = f(p^\alpha)$ for each prime $p (p \leq n)$, $f^p(m) = f(p^\alpha)$ if $p^\alpha \mid m$ and otherwise $f^{(p)}(m) = 0$.

Then $f^p(m) = f(p^\alpha)$ with probability $\frac{1}{n} \left\lfloor \frac{n}{p^\alpha} \right\rfloor$ and $f^{(p)}(m) = 0$ with probability $1 - \frac{1}{n} \left\lfloor \frac{n}{p^\alpha} \right\rfloor$.

Therefore, the average value $f^{(p)}(m)$ over the interval $[1, n]$ is:

$$E[f^p, n] = \frac{f(p^\alpha)}{n} \left\lfloor \frac{n}{p^\alpha} \right\rfloor. \qquad (2.10)$$

It is performed for an additive arithmetic function $f(m), m = 1, ..., n$:

$$f(m) = \sum_{p \leq n} f^{(p)}(m). \qquad (2.11)$$

Based on (2.10), (2.11), the average value of $f(m), m = 1, ..., n$ is:

$$E[f, n] = \sum_{p \leq n} \frac{f(p^\alpha)}{n} \left\lfloor \frac{n}{p^\alpha} \right\rfloor.$$

This implies the desired asymptotics of the mean value of the additive arithmetic function $f(m), m = 1, ..., n$ for the value $n \to \infty$:

$$A_n = \sum_{p \leq n} \frac{f(p^\alpha)}{p^\alpha},$$

which corresponds to (2.9).

An additive function is strongly additive if it satisfies $f^*(p^\alpha) = f^*(p)$, therefore, based on (2.9), the asymptotic of the mean value of a strongly additive arithmetic function $f^*(m), m = 1,...,n$ at the value $n \to \infty$:

$$A_n^* = \sum_{p \leq n} \frac{f^*(p)}{p}. \qquad (2.12)$$

In general, random variables $f^p(m)$ are dependent.

It was shown in [4] that they are asymptotically independent in some cases and, under this assumption, a heuristic formula was found for the asymptotic of the variance of a complex additive arithmetic function for the value $n \to \infty$:

$$D_n = \sum_{p \leq n} \frac{f^2(p^\alpha)}{p^\alpha}. \qquad (2.13)$$

Based on (2.13), the asymptotic of the variance of a complex strongly additive arithmetic function for the value $n \to \infty$ is:

$$D_n^* = \sum_{p \leq n} \frac{f^2(p)}{p}. \qquad (2.14)$$

However, it is shown in the same work that:

$$\frac{1}{n} \sum_{m \leq n} |f^*(m) - A_n| = DD_n^* \text{ and } \frac{1}{n} \sum_{m \leq n} |f^*(m) - A_n| = D^* D_n^*$$

where $D, D^*$ are constants.

Therefore, formulas (2.13) and (2.14) must be used carefully, checking them using (2.2).

Based on (2.6-2.8), an analog of the law of large numbers is fulfilled for additive and strongly additive arithmetic functions.

One can write an expression for the asymptotic of the additive arithmetic function $f(m), m = 1,...,n$, which is satisfied almost everywhere for the value $n \to \infty$:

$$f(n) = A_n + O(b(n)\sqrt{D_n}), \tag{2.15}$$

where $b(n)$ is a slowly growing function, i.e. $b(n) \leq n^{\xi}$, and $\xi > 0$.

Based on (2.15), if $b(n)\sqrt{D_n} = o(A_n)$, then almost everywhere

$$f(n) = A_n(1 + o(1)), \tag{2.16}$$

otherwise

$$f(n) = O(b(n)\sqrt{D_n}) \tag{2.17}$$

at value $n \to \infty$.

Similarly, we can write an expression for the asymptotic of a strongly additive arithmetic function $f^*(m), m = 1,...,n$, which is satisfied almost everywhere for the value $n \to \infty$:

$$f^*(n) = A_n^* + O(b(n)\sqrt{D_n^*}), \tag{2.18}$$

where $b(n)$ is a slowly growing function, i.e. $b(n) \leq n^{\xi}$, and $\xi > 0$.

Based on (2.18), if $b(n)\sqrt{D_n^*} = o(A_n^*)$, then almost everywhere

$$f(n) = A_n^*(1 + o(1)), \tag{2.19}$$

otherwise

$$f^*(n) = O(b(n)\sqrt{D_n^*}) \tag{2.20}$$

at value $n \to \infty$.

However, to determine the asymptotic almost everywhere for an arithmetic function with the value $n \to \infty$ in some cases, it is not necessary to calculate its variance.

Assertion 2

If the arithmetic function $f(m) \geq 0, m = 1,...,n$ and is limited "on average", i.e. $\frac{1}{n}\sum_{m=1}^{n} f(m) \leq C$, then its asymptotic behavior almost everywhere at $n \to \infty$ is equal to:

$$f(n) = O(b(n)), \qquad (2.21)$$

where $b(n)$ is a slowly growing function.

Proof

The arithmetic function satisfies the Markov inequality when the conditions of the assertion in the probability space are met $(\Omega_n, \mathcal{A}_n, \mathbb{P}_n)$ (see the beginning of the chapter):

$$P_n(f(m) \leq b) \geq 1 - \frac{E[f,n]}{b}. \qquad (2.22)$$

Let us put $b = b(n)$ in (2.22), where $b(n)$ is an indefinitely increasing function at value $n \to \infty$ and get:

$$P_n(f(m) \leq b(n)) \geq 1 - \frac{E[f,n]}{b(n)}. \qquad (2.23)$$

The limit $P_n$ exists at value $n \to \infty$ in (2.23):

$$P_n(f(m) \leq b(n)) = 1. \qquad (2.24)$$

Based on (2.24) it follows that the asymptotic almost everywhere for $f(m), m = 1,...,n$ when the value $n \to \infty$ is equal to:

$$f(n) = O(b(n)),$$

which corresponds to (2.21).

It can be seen from the consideration of (2.9), (2.12-2.14) that the probabilistic characteristics of additive arithmetic functions are determined in terms of sums of functions of prime numbers of a certain type. We will consider the definition of these sums in the next chapter of the work.

3. DETERMINATION OF ASYMPTOTICS OF SOME OF FUNCTIONS OF PRIME NUMBERS

Let's consider the sums of functions of prime numbers of the following form:

$$\sum_{\varphi(p) \leq x} g(p), \tag{3.1}$$

where the function of prime number $\varphi(p)$ has an inverse function $\varphi^{-1}(p)$. Then:

$$\sum_{\varphi(p) \leq x} g(p) = \sum_{p \leq \varphi^{-1}(x)} g(p). \tag{3.2}$$

Let us study the case (3.2), when $\varphi(p) = p^\alpha$ and $\alpha > 0$:

$$\sum_{p^\alpha \leq n, \alpha=1,\ldots,k,\ldots} g(p^\alpha) = \sum_{p \leq n} g(p) + \sum_{p \leq n^{1/2}} g(p^2) + \ldots + \sum_{p \leq n^{1/k}} g(p^k) + \ldots \tag{3.3}$$

at $n \to \infty$.

Using (3.3), we define the asymptotic:

$$\sum_{p^\alpha \leq x, \alpha=1,\ldots,k,\ldots} \frac{f(p^\alpha)}{p^\alpha} = \sum_{p \leq n} \frac{f(p)}{p} + \sum_{p \leq n^{1/2}} \frac{f(p^2)}{p^2} + \ldots + \sum_{p \leq n^{1/k}} \frac{f(p^k)}{p^k} + \ldots \tag{3.4}$$

at $n \to \infty$.

We also define the asymptotic:

$$\sum_{p^\alpha \leq x, \alpha=1,\ldots,k,\ldots} \frac{f^2(p^\alpha)}{p^\alpha} = \sum_{p \leq n} \frac{f^2(p)}{p} + \sum_{p \leq n^{1/2}} \frac{f^2(p^2)}{p^2} + \ldots + \sum_{p \leq n^{1/k}} \frac{f^2(p^k)}{p^k} + \ldots \tag{3.5}$$

at $n \to \infty$.

It is known [6] that the asymptotic:

$$\sum_{p \leq x} \frac{\ln p}{p} = \ln x + O(1). \tag{3.6}$$

Considering that the series $\sum_p \frac{\ln p^k}{p^k}$ - converges for $k \geq 2$, then based on (3.4) and (3.6) we obtain the asymptotic:

$$\sum_{p^\alpha \leq n, \alpha=1,\ldots,k,\ldots} \frac{\ln p^\alpha}{p^\alpha} = \ln n + O(1) \tag{3.7}$$

at $n \to \infty$.

Based on [6], the asymptotic of the expression is:

$$\sum_{p \leq n} \frac{\ln^2 p}{p} = \frac{1}{2}\ln^2 n + O(1). \tag{3.8}$$

Then, taking into account that the series $\sum_{p} \frac{\ln^2(p^\alpha)}{p^\alpha}$ at $\alpha \geq 2$ converges, based on (3.5) and (3.8) we get:

$$\sum_{p^\alpha \leq x, \alpha=1,\ldots,k,\ldots} \frac{\ln^2(p^\alpha)}{p^\alpha} = \sum_{p \leq n} \frac{\ln^2(p)}{p} + \sum_{p \leq n^{1/2}} \frac{\ln^2(p^2)}{p^2} + \ldots + \sum_{p \leq n^{1/k}} \frac{\ln^2(p^k)}{p^k} + \ldots = \frac{1}{2}\ln^2 n + O(1) \tag{3.9}$$

at $n \to \infty$.

## 4. DEFINITION OF THE ASYMPTOTIC BEHAVIOR OF ADDITIVE ARITHMETIC FUNCTIONS

Let us single out a class $T$ of additive arithmetic functions for which the additive arithmetic function $f(m), m = 1, 2, \ldots, n$ and the corresponding strongly additive arithmetic function $f^*(m) = \sum_{p \mid m} f(p), m = 1, 2, \ldots, n$, have the same asymptotics of the mean value and variance at $n \to \infty$.

Let us show that the additive arithmetic function $f(m) = \Omega(m)$ (the number of prime divisors of a number, taking into account their multiplicity) and a strongly additive arithmetic function $f^*(m) = \omega(m)$ belong to the class $T$.

The asymptotic of the mean value of an additive arithmetic function $f(m) = \Omega(m), m = 1, \ldots, n$ with the value $n \to \infty$ based on (2.9) is equal to:

$$\sum_{p^\alpha \leq n, \alpha=1,\ldots,k,\ldots} \frac{\Omega(p^\alpha)}{p^\alpha} = \sum_{p \leq n} \frac{1}{p} + \sum_{p^2 \leq n} \frac{2}{p^2} + \ldots = \ln\ln n + O(1). \tag{4.1}$$

The asymptotic of the mean value of a strongly additive arithmetic function $f^*(m) = \omega(m), m = 1, \ldots, n$ at value $n \to \infty$ is equal to:

$$\sum_{p \leq n} \frac{\omega(p)}{p} = \sum_{p \leq n} \frac{1}{p} + O(1) = \ln\ln n + O(1). \tag{4.2}$$

The asymptotic of the variance of an additive arithmetic function $f(m) = \Omega(m), m = 1,...,n$ at value $n \to \infty$ is equal to:

$$\sum_{p^\alpha \leq n, \alpha=1,...,k,...} \frac{\Omega^2(p^\alpha)}{p^\alpha} = \sum_{p \leq n} \frac{1}{p} + \sum_{p^2 \leq n} \frac{4}{p^2} + ... = \ln\ln n + O(1). \tag{4.3}$$

The asymptotic of the variance of a strongly additive arithmetic function $f^*(m) = \omega(m), m = 1,...,n$ at value $n \to \infty$ is:

$$\sum_{p \leq n} \frac{\omega^2(p)}{p} = \sum_{p \leq n} \frac{1}{p} = \ln\ln n + O(1). \tag{4.4}$$

Based on (4.1), (4.2) and (4.3), (4.4), asymptotics of the mean value and variance coincide for arithmetic functions $f(m), f^*(m)$, so the arithmetic function $f(m) = \Omega(m)$ and the strongly additive arithmetic function $f^*(m) = \omega(m)$ belong to the class $T$.

Based on (2.6), almost everywhere when the value $n \to \infty$ asymptotics hold for the indicated functions $f(m), f^*(m)$ $m = 1,2,...,n$:

$$f(n) = f^*(n) = \ln\ln n + O(b(n)\sqrt{\ln\ln n}) = \ln\ln n + O((\ln\ln n)^{1/2+\xi}), \tag{4.5}$$

where $b(n) = (\ln\ln n)^\xi, \xi > 0.$

Turan proved [6] that if an arithmetic function $f^*(m), m = 1,...,n$ is strongly additive and satisfies the condition for all primes $p$: $0 \leq f(p) < c$ and the mean value of this arithmetic function $A_n \to \infty$ at the value $n \to \infty$, then the following form of the analog of the law of large numbers holds for this function:

$$P_n(|f(m) - A_n| \leq b(n)\sqrt{A_n}) \to 1, n \to \infty. \tag{4.6}$$

Formula (4.6) implies the asymptotics for strongly additive arithmetic functions $f^*(m), m = 1,...,n$, which for all primes $p$ satisfies the condition: $0 \leq f(p) < c$ and $A_n \to \infty$ for the value $n \to \infty$:

$$f^*(n) = A_n + O(b(n)\sqrt{A_n}). \tag{4.7}$$

The asymptotic of additive arithmetic function (4.5) is a special case of the asymptotics (4.7).

Assertion 3

If the condition is met for the additive arithmetic function $f(m), m = 1,...,n$ with the value $n \to \infty$:

$$f(n) = O(\ln(n)), \tag{4.8}$$

then the additive arithmetic function $f(m)$ and the corresponding strongly additive arithmetic function $f^*(m) = \sum_{p|m} f(p)$ belong to the class $T$ and asymptotics of these arithmetic functions coincide almost everywhere.

Proof

It was proved in [7] that when condition (4.8) is satisfied, asymptotics of the mean value and variance of an additive arithmetic function $f(m), m = 1,2,...,n$ coincide with the values of the corresponding strongly additive arithmetic function $f^*(m) = \sum_{p|m} f(p), m = 1,2,...,n$ at the value $n \to \infty$; therefore, arithmetic functions $f(m)$ and $f^*(m)$ belong to the class $T$.

Consequently, based on (2.15) and (2.18) asymptotics of the arithmetic functions themselves at the value $n \to \infty$ coincide almost everywhere.

Let's look at an example. It is required to show that the additive arithmetic function $f(m) = \ln \varphi(m)$ for $m = 1,2,...,n$ and the corresponding strongly additive arithmetic function $f^*(m) = \sum_{p|m} \ln \varphi(p)$ for $m = 1,2,...,n$ belong to the class $T$. Find asymptotics of these functions almost everywhere for $n \to \infty$.

It is known that $\varphi(m) = m \prod_{p|m}(1 - 1/p) \leq m$. Therefore, the following trivial estimate is true:

$$\ln \varphi(n) \leq \ln n, \tag{4.9}$$

when the value $n \to \infty$, i.e. the condition of assertion 3 is satisfied.

Therefore, the additive arithmetic function $f(m) = \ln \varphi(m)$ and the corresponding strongly additive arithmetic function $f^*(m) = \sum_{p|m} \ln \varphi(p) = \sum_{p|m} \ln(p-1)$ belong to the class $T$ and have the asymptotic estimate (4.9) almost everywhere when the value $n \to \infty$.

Let us now estimate almost everywhere for these additive and strongly additive arithmetic function in a different way. Since, based on Assertion 3, the probabilistic characteristics for the additive and strongly additive arithmetic functions coincide in this case, we will find probabilistic characteristics only for the strongly additive function $f^*(m) = \sum_{p|m} \ln \varphi(p) = \sum_{p|m} \ln(p-1)$.

Based on (2.12), the asymptotic of the mean value $f^*(m), m = 1,...,n$ is:

$$A_n^* = \sum_{p \le n} \frac{f^*(p)}{p} = \sum_{p \le n} \frac{\ln(p-1)}{p} = \sum_{p \le n} \frac{\ln p(1-1/p)}{p} = \sum_{p \le n} \frac{\ln p}{p} + \sum_{p \le n} \frac{\ln(1-1/p)}{p}. \quad (4.10)$$

Having in mind that the series $\sum_p \frac{\ln(1-1/p)}{p}$ converges and using (3.6) and (4.10) we get:

$$A_n^* = \sum_{p \le n} \frac{f^*(p)}{p} = \ln n + O(1). \quad (4.11)$$

Based on the empirical formula (2.14), we determine the asymptotic of the variance $f^*(m), m = 1,...,n$:

$$D_n^* = \sum_{p \le n} \frac{f^2(p)}{p} = \sum_{p \le n} \frac{\ln^2(p-1)}{p} = \sum_{p \le n} \frac{(\ln p(1-1/p))^2}{p} = \sum_{p \le n} \frac{\ln^2 p}{p} + 2\sum_{p \le n} \frac{\ln p \ln(1-1/p)}{p} + \sum_{p \le n} \frac{\ln^2(1-1/p)}{p}. \quad (4.12)$$

Having in mind that the last two sums in (4.12) are limited and using (3.8) we get:

$$D_n^* = \sum_{p \le n} \frac{f^2(p)}{p} = 0,5\ln^2(n) + O(1). \quad (4.13)$$

Using (4.13) gives the asymptotic almost everywhere for $f^*(m), m = 1,...,n$ as $n \to \infty$:

$$f^*(m) = A_n^* + O(b(n)\sqrt{D_n^*}) = \ln n + O(\ln^{1+\xi} n) = O(\ln^{1+\xi} n).  \quad (4.14)$$

Let us compare (4.14) with the trivial estimate (4.19). Therefore, the empirical formula (2.14) cannot be used in this case.

If we determine the asymptotic of the variance by formula (2.2) in this case, then we will obtain $D_n^* = O(1)$ and an asymptotic estimate almost everywhere for $f^*(m), m = 1,...,n$ when $n \to \infty$:

$$f^*(m) = A_n^* + O(b(n)\sqrt{D_n^*}) = \ln n + O(b(n)),  \quad (4.15)$$

where $b(n)$ is a slowly growing function.

Assertion 4

Let a strongly additive arithmetic function $f^*(m) = \sum_{p|m} f(p) > 0, m = 1,...,n$ and the series $\sum_p \dfrac{f(p)}{p}$ - converge, then the asymptotic almost everywhere $f^*(m), m = 1,...,n$ at value $n \to \infty$ will be equal to:

$$f^*(n) = O(b(n)),  \quad (4.16)$$

where $b(n)$ is a slowly growing function.

Proof

If $n \geq p$, then the number divisible by $p$ among factors $1, 2,..., n$, will be equal to $\left\lfloor \dfrac{n}{p} \right\rfloor$ (Theorem 49) in [8], so the average value $f^*(m), m = 1,...,n$ can be written as:

$$\frac{1}{n}\sum_{p \leq n}\sum_{p|n} f(p) = \frac{1}{n}\sum_{p \leq n} f(p)\left\lfloor \frac{n}{p} \right\rfloor \leq A,  \quad (4.17)$$

where $A = \sum_{p \leq n} \dfrac{f(p)}{p}$.

Based on (4.17), if the series $A = \sum_p \dfrac{f(p)}{p}$ converges, then the average value $f^*(m), m = 1,...,n$ is bounded.

Having in mind that $f^*(m) \geq 0$, then all the conditions of Assertion 2 are satisfied, so the asymptotic is true almost everywhere for the value $n \to \infty$:

$$f^*(n) = O(b(n)),$$

which corresponds to (4.16).

Let's look at an example for assertion 4.

A strongly additive function $f^*(m) = \ln \dfrac{m}{\varphi(m)} \geq 0$, since $m \geq \varphi(m)$.

The series

$$\sum_p \frac{f(p)}{p} = \sum_p \frac{\ln p / \varphi(p)}{p} = \sum_p \frac{\ln p / (p-1)}{p} = \sum_p \frac{-\ln(1-1/p)}{p} = \sum_p \frac{1}{p^2} + o\left(\sum_p \frac{1}{p^2}\right) \text{ converges},$$

therefore, all the conditions of Assertion 4 are satisfied, and the asymptotic almost everywhere $\ln \dfrac{m}{\varphi(m)}$ for $n \to \infty$ is equal to:

$$\ln \frac{m}{\varphi(m)} = O(b(n)),$$

where $b(n)$ is a slowly growing function.

5. ESTIMATE OF ASYMPTOTICS OF MULTIPLIKATIVE ARITHMETIC FUNCTIONS

Let us single out a class $M$ of multiplicative arithmetic functions for which the multiplicative arithmetic function $g(m), m = 1,...,n$ and the corresponding strongly multiplicative arithmetic function $g^*m) = \prod_{p|m} g(p), m = 1,...,n$ have the same asymptotics almost everywhere for the value $n \to \infty$.

Assertion 5

If the condition is met for the multiplicative arithmetic function $g(m), m = 1,...,n$ with the value $n \to \infty$:

$$0 < g(n) \leq n^u, \tag{5.1}$$

where $u > 0$ is a constant, then the multiplicative arithmetic function $g(m), m = 1,...,n$ and the corresponding strongly multiplicative arithmetic function $g^*(m) = \prod_{p|m} g(p), m = 1,...,n$ belong to the class $M$.

Proof

Let's look at the additive arithmetic function $f(m) = \ln(g(m)), m = 1,...,n$.

Based on (5.1), the estimate is made for this arithmetic function: $f(n) = \ln(g(n)) \leq \ln n^u$, i.e. the asymptotic is satisfied $f(n) = O(\ln(n))$ for it, when the value $n \to \infty$.

Based on Assertion 3, if an additive arithmetic function $f(m), m = 1,...,n$ has the asymptotic $f(n) = O(\ln(n))$ at value $n \to \infty$, then this additive function and the corresponding strongly additive function $f^*(m) = \sum_{p|m} f(p)$ have almost everywhere the same asymptotic at the value $n \to \infty$.

Therefore, a strongly multiplicative arithmetic function $g^*(m) = \prod_{p|m} g(p), m = 1,...,n$ has almost everywhere the same asymptotic as a multiplicative function $g(m), m = 1,...,n$ at the value $n \to \infty$ and belongs to the class $M$.

Let's look at examples for assertion 5.

Let it be required to determine the asymptotic almost everywhere of a strongly multiplicative arithmetic function $g^*(m) = \prod_{p|m} p^u(1 - 1/p), m = 1,...,n$ for the value $n \to \infty$.

Для соответствующей мультипликативной арифметической функции $g(m) = m^u \prod_{p|m}(1 - 1/p), m = 1,...,n$ справедлива оценка: $g(n) \leq n^u$, поэтому, на основании утверждения 4, для сильно мультипликативной функции $g^*(m) = \prod_{p|m} p^u(1 - 1/p), m = 1,...,n$ почти всюду при значении $n \to \infty$ справедлива такая же оценка $g^*(n) \leq n^u$.

The estimate is true for the corresponding multiplicative arithmetic function $g(m) = m^u \prod_{p|m}(1 - 1/p), m = 1,...,n$: $g(n) \leq n^u$, therefore, based on Assertion 4, the same

estimate is true almost everywhere with the value $n \to \infty$ for a strongly multiplicative function

$$g^*(m) = \prod_{p|m} p^u(1-1/p), m = 1,...,n - g^*(n) \leq n^u. \quad (5.2)$$

Let's consider another example. Let it be required to find the asymptotic of a strongly multiplicative arithmetic function $g^*(m) = e^{\omega(m)+\sum_{p|m}\ln(1-1/p)}, m = 1,...,n, n \to \infty$.

Here, the corresponding strongly additive arithmetic function $f^*(m) = \omega(m) + \sum_{p|m}\ln(1-1/p), m = 1,...,n, n \to \infty$ satisfies the conditions of Turan's assertion, since $\omega(p) = 1$. Therefore, the asymptotic $f^*(m)$ can be determined from formula (4.7).

It suffices to define the asymptotic of the mean value $f^*(m), m = 1,...,n, n \to \infty$:

$$A_n = \sum_{p \leq n} \frac{\omega(p)}{p}(1 - \frac{1}{p}) + \sum_{p \leq n} \frac{\ln(1-1/p)}{p}(1 - \frac{1}{p}) = \sum_{p \leq n}\frac{1}{p} - \sum_{p \leq n}\frac{1}{p^2} + O(1) = \ln\ln n + O(1). \quad (5.3)$$

Based on (4.7), (5.3), the asymptotic almost everywhere for $f^*(m), m = 1,...,n, n \to \infty$ holds:

$$f^*(n) = \ln\ln n + O(b(n)\sqrt{\ln\ln n}) = \ln\ln n + O((\ln\ln n)^{1/2+\xi}), \quad (5.4)$$

where $b(n) = (\ln\ln n)^\xi, \xi > 0$.

Having in mind (5.4), we can obtain an almost everywhere estimate for the corresponding strongly multiplicative function $g^*(m) = e^{\omega(m)+\sum_{p|m}\ln(1-1/p)}, m = 1,...,n, n \to \infty$:

$$g^*(n) = e^{\ln\ln n + O((\ln\ln n)^{1/2+\xi})} = \ln(n)e^{O((\ln\ln n)^{1/2+\xi})}. \quad (5.5)$$

Assertion 6

Let $g^*(m) = \prod_{p|m} g(p) \geq 1, m = 1,...,n$ is a strongly multiplicative arithmetic function and also the series $\sum_p \ln g^*(p)$ converge, then the asymptotic almost everywhere $g^*(m), m = 1,...,n$ at the value $n \to \infty$ is equal to:

$$g^*(m) = \prod_{p|m} g(p) = O(b(n)), \quad (5.6)$$

where $b(n)$ is a slowly growing function.

Proof

Let's consider a strongly additive arithmetic function $\ln g^*(m) = \ln \prod_{p|m} g(p) = \sum_{p|m} \ln g(p)$.

Having in mind that $g^*(m) \geq 1$, then $\ln g^*(m) \geq 0$.

Since the condition that the series $\sum_{p} \ln g^*(p)$ converges is satisfied, then all the conditions of Assertion 4 are satisfied for a strongly additive arithmetic function $\ln g^*(m), m = 1,...,n$, and therefore the asymptotic almost everywhere for the value $n \to \infty$ is equal to:

$$\ln g^*(m) = \sum_{p|m} \ln g(p) = O(b_1(n)),  \qquad (5.7)$$

where $b_1(n)$ is a slowly growing function.

Based on (5.7), the asymptotic almost everywhere for the value $n \to \infty$ is equal to:

$$e^{\ln g^*(m)} = g^*(m) = e^{O(b_1(n))} = O(b(n)), \qquad (5.8)$$

where $b(n)$ is a slowly growing function.

Estimate (5.8) corresponds to (5.6).

Let's look at an example for assertion 6.

$g^*(m) = \dfrac{m}{\varphi(m)} \geq 1$ is the strongly multiplicative arithmetic function.

On the other hand, the series $\sum_{p|m} \ln \dfrac{p}{\varphi(p)}$ - converges (see the example for Assertion 4).

Therefore, the conditions of Assertion 6 are satisfied for $g^*(m) = \dfrac{m}{\varphi(m)}$.

Therefore, the asymptotic is true almost everywhere, when the value $n \to \infty$, for this strongly multiplicative arithmetic function:

$$g^*(n) = \frac{n}{\varphi(n)} = O(b(n)),$$

where $b(n)$ is a slowly growing function.

### 6. CONCLUSION AND SUGGESTIONS FOR FURTHER WORK

The next article will continue to study the asymptotic behavior of some arithmetic functions.

### 7. ACKNOWLEDGEMENTS

Thanks to everyone who has contributed to the discussion of this paper. I am grateful to everyone who expressed their suggestions and comments in the course of this work.

# References

1. J.Chattopadhyay, P. Darbar Mean values and moments of arithmetic functions over number fields, Research in Number Theory 5(3), 2019.

2. M. Garaev, M. Kühleitner, F. Luca, G. Nowak Asymptotic formulas for certain arithmetic functions, Mathematica Slovaca 58(3):301-308, 2008.

3. R. Brad Arithmetic functions in short intervals and the symmetric group, arXiv preprint https://arxiv.org/abs/1609.02967(2016).

4. J. Kubilius. Probabilistic methods in number theory. Vilnius, 1962, 220 p.

5. Volfson V.L. Asymptotics of sums of functions of primes located on an arithmetic progression, arXiv preprint https://arxiv.org/abs/2112.04474 (2021)

6. P.Turan. Uber einige Verungen eines Satzes von Hardy und Ramanujan. J.London Math. Soc., 1936, 11, 125-133.

7. Volfson V.L. Asymptotics of probability characteristics of additive arithmetic functions, arXiv preprint https://arxiv.org/abs/2108.13199 (2021)

8. A.A. Bukhshtab. "Theory of numbers", Publishing house "Education", M., 1966, 384 p.